\documentclass{article}

\usepackage[top=2cm, bottom=2cm, left=3cm, right=3cm]{geometry}

\usepackage{amsmath, amsthm, bm, mathrsfs, mathtools, amssymb, amsfonts}
\usepackage{graphicx}
\usepackage[colorlinks=true, allcolors=blue]{hyperref}
\usepackage{multicol}
\usepackage[all]{xy}
\usepackage{appendix}

\newtheorem{theorem}{Theorem}[section]
\newtheorem{lemma}[theorem]{Lemma}
\newtheorem{proposition}[theorem]{Proposition}

\theoremstyle{definition}

\theoremstyle{remark}
\newtheorem{remark}{Remark}[section]

\newcommand{\Cov}{\operatorname{Cov}}
\newcommand{\Var}{\operatorname{Var}}

\allowdisplaybreaks

\begin{document}

\title{The maximal correlation coefficient associated with the minimum}
\author{Yinshan Chang\thanks{Address: College of Mathematics, Sichuan University, Chengdu 610065, China; Email: ychang@scu.edu.cn.}, Qinwei Chen\thanks{Address: College of Mathematics, Sichuan University, Chengdu 610065, China; Email: qinweic@outlook.com.}}

\date{}
\maketitle

\begin{abstract}
For independent random variables $(X_i)_{1\leq i\leq n}$, we consider the maximal correlation coefficient $R=R(\min_{i:1\leq i\leq m}X_i,\min_{j:\ell+1\leq j\leq n}X_j)$. If $X_1,X_2,\ldots,X_n$ are identically distributed with the same continuous distribution, we find that $R=(m-\ell)/\sqrt{m(n-\ell)}$. For discrete distributions, we calculate the maximal correlation coefficient $R$ for Bernoulli distributions, geometric distributions, binomial distributions and Poisson distributions. Our paper answers a question in \cite[Section~6]{ChangChen}.
\end{abstract}

\section{Introduction}
The maximal correlation coefficient introduced by Gebelein in \cite{GebeleinMR7220} measures the dependency of two random variables. For two non-degenerate random variables $X$ and $Y$, it is given by
\begin{equation*}
R(X,Y)=\sup\frac{\Cov(\varphi(X),\psi(Y))}{\sqrt{\Var(\varphi(X))\Var(\psi(Y))}},
\end{equation*}
where the supremum is taken over all measurable functions $\varphi$ and $\psi$ such that $\Var(\varphi(X)),\Var(\psi(Y))\in(0,\infty)$.

For the literature on the maximal correlation coefficient, we refer to our recent paper \cite{ChangChen} and the references therein.

Among various results, B\"{u}cher and Staud \cite[Corollary~2.2]{BucherStaudMR4835998} have calculated the maximal correlation coefficient for the bivariate Marshall-Olkin exponential distribution. To be more precise, consider independent exponential random variables $X_1,X_2$ and $X_3$ with parameters $\lambda_1,\lambda_2$ and $\lambda_3$, respectively. Then
\begin{equation*}
R(\min(X_1,X_2),\min(X_2,X_3))=\frac{\lambda_2}{\sqrt{(\lambda_1+\lambda_2)(\lambda_2+\lambda_3)}}.
\end{equation*}
In \cite[Section~5.6]{ChangChen}, we reproduce the upper bound
\begin{equation}\label{eq: upper bound bivariate M-O}
R(\min(X_1,X_2),\min(X_2,X_3))\leq\frac{\lambda_2}{\sqrt{(\lambda_1+\lambda_2)(\lambda_2+\lambda_3)}}
\end{equation}
of their work, which follows as a consequence of the following proposition.
\begin{proposition}[Proposition~5.5 in \cite{ChangChen}]\label{prop: Yu replace sum by min}
Let $X_1,X_2,\ldots,X_n$ be i.i.d. (real) random variables. Let $1\leq \ell+1\leq m\leq n$. Then, we have that
\begin{equation}\label{eq: Yu replace sum by min}
R\left(\min_{i:1\leq i\leq m}X_i,\min_{j:\ell+1\leq j\leq n}X_j\right)\leq\frac{m-\ell}{\sqrt{m(n-\ell)}}.
\end{equation}
\end{proposition}
When each $X_i$ follows an exponential distribution with parameter $\lambda$, we find that
\begin{align*}
&R\left(\min_{i:1\leq i\leq m} X_i,\min_{j:\ell+1\leq j\leq n} X_j\right)
\\&= R\left(\min\left\{\min_{i:1\leq i\leq \ell} X_i,\min_{i:\ell+1\leq i\leq m} X_i\right\},
\min\left\{\min_{j:\ell+1\leq j\leq m} X_j, \min_{j:m+1\leq j\leq n} X_j\right\}\right)\\
&= R\left(\min(V_1, V_2),\min(V_2, V_3)\right),
\end{align*}
where $V_1, V_2, V_3$ are independent exponential random variables with parameters $\lambda  \ell$, $\lambda(m-\ell)$, and $\lambda(n-m)$, respectively. Then by \cite[Corollary~2.2]{BucherStaudMR4835998}, the upper bounds in \eqref{eq: upper bound bivariate M-O} and \eqref{eq: Yu replace sum by min} are sharp when each $X_i$ follows an exponential distribution. We thus wonder whether they are also sharp for other distributions, see the second open problem in \cite[Section~6]{ChangChen}.

In the present paper, we prove that \eqref{eq: Yu replace sum by min} is also sharp for other continuous distributions. However, it is not sharp for some discrete distributions. Moreover, we obtain explicit formulae for Bernoulli, geometric, binomial and Poisson distributions.
\begin{theorem}\label{thm: continuous random variables}
Let $X_1,X_2,\ldots,X_n$ be i.i.d. (real) continuous random variables. Let $1\leq \ell+1\leq m\leq n$. Then we have that
\begin{equation*}
R\left(\min_{i:1\leq i\leq m}X_i,\min_{j:\ell+1\leq j\leq n}X_j\right)=\frac{m-\ell}{\sqrt{m(n-\ell)}}.
\end{equation*}
\end{theorem}

\begin{theorem}\label{thm: Bernoulli}
Let $X_1,X_2,\ldots,X_n$ be independent Bernoulli random variables such that $P(X_i=1)=1-P(X_i=0)=p_i\in(0,1)$ for $1\leq i\leq n$. Let $R=R\left(\min\limits_{i:1\leq i\leq m}X_i,\min\limits_{j:\ell+1\leq j\leq n}X_j\right)$, where $1\leq \ell+1\leq m\leq n$. Then
\begin{equation}\label{eq:Bernoulli}
R=\frac{\prod_{i=1}^{n}p_i-\left(\prod_{i=1}^{m}p_i\right)\left(\prod_{i=\ell+1}^{n}p_i\right)}{\sqrt{\left(\prod_{i=1}^{m}p_i\right)\left(1-\prod_{i=1}^{m}p_i\right)\left(\prod_{i=\ell+1}^{n}p_i\right)\left(1-\prod_{i=\ell+1}^{n}p_i\right)}}.
\end{equation}
In particular, if $p_1=p_2=\cdots=p_n=p\in (0,1)$, we have $R=R_{m,\ell}(p)$, where
\begin{equation}\label{eq: R(p)}
R_{m,\ell}(p):=\frac{p^n-p^{m+n-\ell}}{\sqrt{p^{m+n-\ell}(1-p^m)(1-p^{n-\ell})}}.
\end{equation}
Moreover, $R_{m,\ell}(p)$ is non-decreasing in $p$.
\end{theorem}

\begin{theorem}\label{thm: Geometric}
Let $X_1,X_2,\ldots,X_n$ be independent geometric random variables. For $1\leq i\leq n$, the parameter of $X_i$ is $p_i\in(0,1)$. Let $1\leq \ell+1\leq m\leq n$. Then we have that
\begin{align*}
&R\left(\min_{i:1\leq i\leq m}X_i,\min_{j:\ell+1\leq j\leq n}X_j\right)=\\ &\frac{\prod_{i=1}^{n}\left(1-p_i\right)-\left(\prod_{i=1}^{m}\left(1-p_i\right)\right)\left(\prod_{i=\ell+1}^{n}\left(1-p_i\right)\right)}{\sqrt{\left(\prod_{i=1}^{m}\left(1-p_i\right)\right)\left(1-\prod_{i=1}^{m}\left(1-p_i\right)\right)\left(\prod_{i=\ell+1}^{n}\left(1-p_i\right)\right)\left(1-\prod_{i=\ell+1}^{n}\left(1-p_i\right)\right)}}.
\end{align*}
In particular, for $n=3$, $m=2$ and $\ell=1$, we have that
\begin{equation*}
\begin{aligned}
&R(\min(X_1,X_2),\min(X_2,X_3))\\
&=\frac{\sqrt{(1-p_1)(1-p_3)}p_2}{\sqrt{(1-(1-p_1)(1-p_2))(1-(1-p_2)(1-p_3))}}. \end{aligned}
\end{equation*}
\end{theorem}

\begin{theorem}\label{thm: binomial}
Let $X_1,X_2,\ldots,X_n$ be i.i.d. binomial random variables with parameter $(d,p)$ with $0<p<1$. Let $1\leq \ell+1\leq m\leq n$. Then we have that
\begin{equation}
R\left(\min_{i:1\leq i\leq m}X_i,\min_{j:\ell+1\leq j\leq n}X_j\right)=R_{m,\ell}(1-(1-p)^d),
\end{equation}
where $R_{m,\ell}(p)$ is defined in \eqref{eq: R(p)}.
\end{theorem}

\begin{theorem}\label{thm: Poisson}
Let $X_1,X_2,\ldots,X_n$ be i.i.d. Poisson random variables with parameter $\lambda>0$. Let $1\leq \ell+1\leq m\leq n$. Then we have that
\begin{equation}\label{eq: Poisson}
R\left(\min_{i:1\leq i\leq m}X_i,\min_{j:\ell+1\leq j\leq n}X_j\right)=R_{m,\ell}(1-e^{-\lambda}),
\end{equation}
where $R_{m,\ell}(p)$ is defined in \eqref{eq: R(p)}.
\end{theorem}

\emph{Organization of the paper}: In Section~\ref{sect: preliminaries}, we present three useful properties of the maximal correlation coefficient. In Section~\ref{sect: Continuous distribution}, we prove Theorem~\ref{thm: continuous random variables}. In Section~\ref{sect:special discrete distributions}, we first prove Theorem~\ref{thm: Bernoulli}. Then we prove Theorems~\ref{thm: Geometric} and \ref{thm: binomial} by using Theorem~\ref{thm: Bernoulli} and prove Theorem~\ref{thm: Poisson} by using Theorems~\ref{thm: Bernoulli} and \ref{thm: binomial}.

\section{Preliminary}\label{sect: preliminaries}
We present three useful properties of the maximal correlation coefficient.

If $U=\varphi(X)$ is a measurable function of $X$ and $V=\psi(Y)$ is a measurable function of $Y$, then $R(U,V)\leq R(X,Y)$.

Suppose that $(X_k,Y_k)$ converges to $(X,Y)$ in law in some proper space. By \cite[Lemma~2.3]{ChangChen},
\begin{equation}\label{eq: lower semicontinuity}
R(X,Y)\leq\liminf_{k\to\infty}R(X_k,Y_k).
\end{equation}

\begin{lemma}\label{lem: Csaki-Fischer}
Consider a sequence of independent random variables $(X_i,Y_i)_{i\geq 1}$, and suppose that $X_i$ takes values in some Polish space $\mathcal{X}_i$, while $Y_i$ takes values in some Polish space $\mathcal{Y}_i$. Then we have that
\[R((X_i)_{i\geq 1},(Y_i)_{i\geq 1})=\sup_{i\geq 1}R(X_i,Y_i).\]
\end{lemma}
\begin{proof}
Suppose that $d_{X,i}(\cdot,\cdot)$ is a metric on the Polish space $\mathcal{X}_i$, and $d_{Y,i}(\cdot,\cdot)$ is a metric on the Polish space $\mathcal{Y}_i$. Then we define a metric $d_{X}(\cdot,\cdot)$ on the product space  $\mathcal{X}=\prod_{i\geq 1}\mathcal{X}_i$: For $x=(x_i)_{i\geq 1}$ and $\tilde{x}=(\tilde{x}_i)_{i\geq 1}$, we define
\[d_{X}(x,\tilde{x})=\sum_{i=1}^{\infty}\min(d_{X,i}(x_i,\tilde{x}_i),1)/2^{i}.\]
Then $(\mathcal{X},d_{X})$ is a product Polish space. Similarly, we construct the product Polish space $(\mathcal{Y},d_{Y})$. Then for each $i\geq 1$, we fix a particular point $x_i^{0}$ in the Polish space $\mathcal{X}_i$ and a particular point $y_i^{0}$ in the Polish space $\mathcal{Y}_i$. Let $((X_i)_{i\geq 1},(Y_i)_{i\geq 1})$ have the joint distribution $\mu$ on the product Polish space $\mathcal{X}\times\mathcal{Y}$. For $m\geq 1$, define \[\boldsymbol{X}^{(m)}=(X_1,X_2,\cdots,X_m,x_{m+1}^{0},x_{m+2}^{0},\ldots)\text{ and }\boldsymbol{Y}^{(m)}=(Y_1,Y_2,\cdots,Y_m,y_{m+1}^{0},y_{m+2}^{0},\ldots),\]
with the joint distribution $\mu_m$. Then by Cs\'{a}ki-Fischer identity (see \cite[Theorem~6.2]{CsakiFischerMR0166833} or \cite[Theorem~1]{WitsenhausenMR363678}),
\[R(\mu_m)=R\left(\boldsymbol{X}^{(m)},\boldsymbol{Y}^{(m)}\right)=R\left((X_i)_{1\leq i\leq m},(Y_i)_{1\leq i\leq m}\right)=\max_{1\leq i\leq m}R(X_i,Y_i).\]
As $m\to\infty$, $\left(\boldsymbol{X}^{(m)},\boldsymbol{Y}^{(m)}\right)$ converges almost surely to $((X_i)_{i\geq 1},(Y_i)_{i\geq 1})$. Hence, as $m\to\infty$, the joint distribution $\mu_m$ of $\left(\boldsymbol{X}^{(m)},\boldsymbol{Y}^{(m)}\right)$ converges weakly to the joint distribution $\mu$ of $((X_i)_{i\geq 1},(Y_i)_{i\geq 1})$. By \eqref{eq: lower semicontinuity} we have that,
\[
R\left((X_i)_{i\geq 1},(Y_i)_{i\geq 1}\right)\leq\liminf_{m\to\infty}R\left(\boldsymbol{X}^{(m)},\boldsymbol{Y}^{(m)}\right)=\sup_{i\geq 1}R(X_i,Y_i).
\]
On the other hand, for each $i\geq 1$, $X_i$ is a measurable function of $(X_j)_{j\geq 1}$ and $Y_i$ is a measurable function of $(Y_j)_{j\geq 1}$. Hence, we have that
\[R(X_i,Y_i)\leq R\left((X_j)_{j\geq 1},(Y_j)_{j\geq 1}\right).\]
Then we have the conclusion.
\end{proof}

\section{Continuous distribution}\label{sect: Continuous distribution}
As explained in the introduction, the upper bound in \eqref{eq: Yu replace sum by min} is sharp when each $X_i$ follows an exponential distribution. For general continuous random variables $X_1,X_2,\ldots,X_n$ with the common distribution function $F(x)$, define $Y_i=\varphi(X_i)$, where $\varphi(x)=-\ln(1-F(x))$. Note that $\varphi$ is a non-decreasing function. Moreover, $Y_1,Y_2,\ldots,Y_n$ are i.i.d. exponential random variables with parameter $1$. Hence, we have that
\begin{align}\label{eq: lower bound}
R\left(\min_{i:1\leq i\leq m}X_i,\min_{j:\ell+1\leq j\leq n}X_j\right)&\geq R\left(\varphi\left(\min_{i:1\leq i\leq m}X_i\right),\varphi\left(\min_{j:\ell+1\leq j\leq n}X_j\right)\right)\notag\\
& = R\left(\min_{i:1\leq i\leq m}\varphi(X_i),\min_{j:\ell+1\leq j\leq n}\varphi(X_j)\right)\notag\\
& = R\left(\min_{i:1\leq i\leq m}Y_i,\min_{j:\ell+1\leq j\leq n}Y_j\right)\notag\\
& = \frac{m-\ell}{\sqrt{m(n-\ell)}}.
\end{align}
By \eqref{eq: Yu replace sum by min} and \eqref{eq: lower bound}, the proof is complete.
\begin{remark}
Alternatively, the lower bound \eqref{eq: lower bound} could be derived from Theorem~\ref{thm: Bernoulli}. Indeed, let $h(x)=1_{(\delta(p),\infty)}(x)$, where $P(X_i>\delta(p))=p$. Then
\begin{align*}
R\left(\min_{i:1\leq i\leq m}X_i,\min_{j:\ell+1\leq j\leq n}X_j\right)&\geq R\left(h\left(\min_{i:1\leq i\leq m}X_i\right),h\left(\min_{j:\ell+1\leq j\leq n}X_j\right)\right)\notag\\
& = R\left(\min_{i:1\leq i\leq m}h(X_i),\min_{j:\ell+1\leq j\leq n}h(X_j)\right)\notag\\
& \overset{Theorem~\ref{thm: Bernoulli}} = R_{m,\ell}(p)\\
&\overset{p\to 1}{\to}\frac{m-\ell}{\sqrt{m(n-\ell)}}.
\end{align*}
\end{remark}

\section{Some special discrete distributions}\label{sect:special discrete distributions}
\subsection{Bernoulli distribution}
In this subsection, we prove Theorem~\ref{thm: Bernoulli}.

Firstly, we calculate the maximal correlation coefficient between two Bernoulli random variables. The following lemma is taken from \cite[Section~3.2]{ChangChen}.
\begin{lemma}\label{lemma 1}
Consider a random vector $\left(\xi,\eta\right)$ taking values in $\{0,1\}^2$, with the joint distribution in Table~\ref{table: joint distribution}.
\begin{table}
\centering
\begin{tabular}{|c|c|c|}
\hline
$\xi\backslash\eta$ & $0$ & $1$ \\
\hline
$0$ & $p_{00}$ & $p_{01}$ \\
\hline
$1$ & $p_{10}$ & $p_{11}$ \\
\hline
\end{tabular}
\caption{The joint distribution of $(\xi,\eta)$}
\label{table: joint distribution}
\end{table}
Then we have that
\begin{equation*}
R(\xi,\eta)=\left.\left|\det\begin{pmatrix}p_{00}&p_{01}\\p_{10}&p_{11}\end{pmatrix}\right.\right|/\sqrt{(p_{00}+p_{01})(p_{10}+p_{11})(p_{00}+p_{10})(p_{01}+p_{11})}.
\end{equation*}
\end{lemma}

Next, consider independent Bernoulli random variables $X_1,X_2,\ldots,X_n$, such that $P(X_i=1)=1-P(X_i=0)=p_i$ for $1\leq i\leq n$. Let $X=\min_{i:1\leq i\leq m}X_i$, $Y=\min_{j:\ell+1\leq j\leq n}X_j$. Then $X$ and $Y$ are both Bernoulli random variables with the joint distribution described in Table~\ref{table: Bernoulli joint distribution}.

\begin{table}
\centering
\begin{tabular}{|c|c|c|}
\hline
$X\backslash Y$ & $0$ & $1$ \\
\hline
$0$ & $ 1-{\textstyle \prod_{i=\ell+1}^{n}}p_i-{\textstyle \prod_{i=1}^{m}}p_i +{\textstyle \prod_{i=1}^{n}}p_i$ & ${\textstyle \prod_{i=\ell+1}^{n}}p_i- {\textstyle \prod_{i=1}^{n}}p_i $ \\
\hline
$1$ & ${\textstyle \prod_{i=1}^{m}}p_i- {\textstyle \prod_{i=1}^{n}}p_i$ & $ {\textstyle\prod_{i=1}^{n}}p_i$ \\
\hline
\end{tabular}
\caption{The joint distribution of $(X,Y)$}
\label{table: Bernoulli joint distribution}
\end{table}
By Lemma~\ref{lemma 1}, we have the \eqref{eq:Bernoulli}. Finally, we prove the the monotonicity of $p\mapsto R_{m,\ell}(p)$. let $a=(n+\ell-m)/2$, $b=m-\ell$ and $c=n-\ell$. Then for $0<p<1$,
\begin{align*}
\frac{R_{m,\ell}'(p)}{R_{m,\ell}(p)} &= \frac{1}{p}\left(a-\frac{bp^{b}}{1-p^{b}} + \frac{m p^{m}}{2(1 - p^m)} + \frac{cp^{c}}{2(1 - p^{c})}\right)\\
&\overset{2a+2b-m-c=0}{=} \frac{1}{2p}\left(\left(\frac{m}{1-p^{m}}-\frac{b}{1-p^b}\right)+\left(\frac{c}{1-p^{c}}-\frac{b}{1-p^b}\right)\right).
\end{align*}
For $x>0$, define $f(x)=x/(1-p^{x})$. Then $f'(x)=\frac{1-p^{x}+xp^{x}\ln p}{(1-p^{x})^2}>0$.
Hence, $f(x)$ is increasing. Since $b\leq m$, we have that $f(m)-f(b)\geq 0$. Similarly, $f(c)-f(b)\geq 0$. Therefore, we have that
\begin{equation*}
\frac{R_{m,\ell}'(p)}{R_{m,\ell}(p)}=\frac{1}{2p}(f(m)-f(b)+f(c)-f(b))\geq 0,
\end{equation*}
and $R_{m,\ell}(p)$ is non-decreasing in $p$.

\subsection{Geometric distribution}
In this subsection, we prove Theorem~\ref{thm: Geometric}. Since the geometric distribution and the Bernoulli distribution are closely related, we will deduce Theorem~\ref{thm: Geometric} from Theorem~\ref{thm: Bernoulli}.

Consider independent Bernoulli random variables $(B_{ik})_{i,k\geq 1}$ such that $B_{ik}$ has the parameter $1-p_i$ for $i,k\geq 1$. For $i\geq 1$, let $X_i=\min\{k\geq 1:B_{ik}=0\}$. Then $(X_i)_{i\geq 1}$ is a sequence of independent geometric random variables such that $X_i$ has the parameter $p_i$ for $i\geq 1$.

Let $U=\min_{i:1\leq i\leq m}X_i$, $V=\min_{j:\ell+1\leq j\leq n}X_j$, and $U_k=\min_{i:1\leq i\leq m}B_{ik}$, we have that $U=\min\{k\geq 1:U_{k}=0\}$. Similarly, let $V_k=\min_{j:\ell+1\leq j\leq n}B_{jk}$, we have that $V=\min\{k\geq 1:V_{k}=0\}$. On the other hand, $U_{1}=1_{U\geq 2}\text{ and }V_{1}=1_{V\geq 2}$.

By Theorem~\ref{thm: Bernoulli}, for $k\geq 1$, we have that
\begin{align*}
&R(U_k,V_k)=R(U_1,V_1)=\\&\frac{\prod_{i=1}^{n}\left(1-p_i\right)-\left(\prod_{i=1}^{m}\left(1-p_i\right)\right)\left(\prod_{i=\ell+1}^{n}(1-p_i)\right)}{\sqrt{\left(\prod_{i=1}^{m}\left(1-p_i\right)\right)\left(1-\prod_{i=1}^{m}\left(1-p_i\right)\right)\left(\prod_{i=\ell+1}^{n}(1-p_i)\right)\left(1-\prod_{i=\ell+1}^{n}\left(1-p_i\right)\right)}}.
\end{align*}

Since $(U_k,V_k)_{k\geq 1}$ is a sequence of i.i.d. random vectors, by Lemma~\ref{lem: Csaki-Fischer}, we have that 
\[R((U_k)_{k\geq 1},(V_k)_{k\geq 1})=\sup_{k\geq 1}R(U_k,V_k)=R(U_1,V_1).\]
Finally, since $U$ is a measurable function of $(U_k)_{k\geq 1}$ and $V$ is a measurable function of $(V_k)_{k\geq 1}$, we have that $R(U,V)\leq R((U_k)_{k\geq 1},(V_k)_{k\geq 1})$. On the other hand, since $U_1$ is a measurable function of $U$ and $V_1$ is a measurable function of $V$, we have that $R(U,V)\geq R(U_1,V_1)$. Finally, we have that $R(U,V)=R(U_1,V_1)$ and the proof is complete.

\subsection{Binomial distribution}
In this subsection, we prove Theorem~\ref{thm: binomial} by using Theorem~\ref{thm: Bernoulli}.

Consider i.i.d. binomial random variables $X_1,X_2,\ldots,X_n$ with parameter $(d,p)$. Let $(B_{ik})_{i,k\geq 1}$ be independent Bernoulli random variables such that
\begin{equation*}
P(B_{ik}=1)=p(k):=P(X_i\geq k|X_i\geq k-1).
\end{equation*}
For $k\geq d+1$ and $i\geq 1$, $B_{ik}=0$.
Then we need the following lemma on the monotonicity of $k\mapsto p(k)$.
\begin{lemma}\label{lem: monotonicity of p(k)}
The conditional probability $p(k)$ is non-increasing in $k$. In particular, for all $k\geq 1$, $p(k)\leq p(1)=1-(1-p)^d$.
\end{lemma}
\begin{proof}
Note that
\begin{equation*}
\frac{p(k+1)}{p(k)}=\frac{P(X_i\geq k+1)P(X_i\geq k-1)}{\left(P(X_i\geq k)\right)^2}.
\end{equation*}
Hence, it suffices to prove that
\begin{equation*}
\Delta:=\left(P(X_i\geq k)\right)^2-P(X_i\geq k+1)P(X_i\geq k-1)\geq 0\text{ for }k\geq 1.
\end{equation*}
Note that
\[\Delta=P(X_i=k)P(X_i\geq k)-P(X_i=k-1)P(X_i\geq k+1).\]
For $k\geq 1$ and $j\geq 0$, since
\begin{equation*}
\binom{n}{k}=\frac{n-k+1}{k}\binom{n}{k-1}\text{ and }\binom{n}{k+j+1}=\frac{n-k-j}{k+j+1}\binom{n}{k+j},
\end{equation*}
we have that
\begin{equation*}
\binom{n}{k}\binom{n}{k+j}\geq\binom{n}{k-1}\binom{n}{k+j+1}.
\end{equation*}
Therefore, we have that
\begin{equation*}
P(X_i=k)P(X_i=k+j)\geq P(X_i=k-1)P(X_i=k+j+1).
\end{equation*}
By summing over $j\geq 0$, we see that $\Delta\geq 0$, and the proof is complete.
\end{proof}
Let us continue with the proof of Theorem~\ref{thm: binomial}. We have the following equality in distribution: $X_i\overset{(d)}{=}\min\{k\geq 1:B_{ik}=0\}-1$.
Without loss of generality, we assume that the above equalities hold almost surely. Let $U=\min_{i:1\leq i\leq m}X_i$, $V=\min_{j:\ell+1\leq j\leq n}X_j$, and let $U_k=\min_{i:1\leq i\leq m}B_{ik}$, we have that $U=\min\{k\geq 1:U_{k}=0\}-1$. Similarly, let $V_k=\min_{j:\ell+1\leq j\leq n}B_{jk}$, we have that $V=\min\{k\geq 1:V_{k}=0\}-1$. On the other hand,
$U_{1}=1_{U>0}\text{ and }V_{1}=1_{V>0}$.

Since $U_1$ is a measurable function of $U$ and $V_1$ is a measurable function of $V$, we have that $R(U_1,V_1)\leq R(U,V)$. Since $U$ is a measurable function of $(U_k)_{k\geq 1}$ and $V$ is a measurable function of $(V_k)_{k\geq 1}$, we have that $R(U,V)\leq R((U_k)_{k\geq 1},(V_k)_{k\geq 1})$. By Lemma~\ref{lem: Csaki-Fischer}, Theorem~\ref{thm: Bernoulli}, and Lemma~\ref{lem: monotonicity of p(k)},
\begin{equation*}
R((U_k)_{k\geq 1},(V_k)_{k\geq 1})=\sup_{k\geq 1}R(U_k,V_k)=\sup_{k\geq 1}R_{m,\ell}(p(k))=R_{m,\ell}(1-(1-p)^d).
\end{equation*}

\subsection{Poisson distribution}
In this subsection, we deduce Theorem~\ref{thm: Poisson} from Theorems~\ref{thm: Bernoulli} and \ref{thm: binomial}.

For $k\geq 1$, let $X^{(k)}_1,X^{(k)}_2,\ldots,X^{(k)}_n$ be i.i.d. binomial random variables with parameter $(k,\lambda/k)$. Then $(X^{(k)}_1,X^{(k)}_2,\ldots,X^{(k)}_n)$ converges in distribution towards $(X_1,X_2,\ldots,X_n)$ as $k\to\infty$, where $X_1,X_2,\ldots,X_n$ are i.i.d. Poisson random variables with parameter $\lambda$. By Theorem~\ref{thm: binomial}, we have that
\begin{equation*}
R\left(\min_{i:1\leq i\leq m}X^{(k)}_i,\min_{j:\ell+1\leq j\leq n}X^{(k)}_j\right)=R_{m,\ell}(1-(1-\lambda/k)^{k}).
\end{equation*}
By taking $k\to\infty$, using \eqref{eq: lower semicontinuity} and the continuity of $R_{m,\ell}$, we obtain that
\begin{equation*}
R\left(\min_{i:1\leq i\leq m}X_i,\min_{j:\ell+1\leq j\leq n}X_j\right)\leq R_{m,\ell}(1-e^{-\lambda}).
\end{equation*}
On the other hand, we take $Y_i=1_{X_i>0}$. Then $Y_1,Y_2,\ldots,Y_n$ are i.i.d. Bernoulli random variables with parameter $1-e^{-\lambda}$. Moreover, the Bernoulli random variable $\min_{i:1\leq i\leq m}Y_i$ is equal to $0$ if and only if $\min_{i:1\leq i\leq m}X_i=0$. Hence, $\min_{i:1\leq i\leq m}Y_i$ is a measurable function of
$\min_{i:1\leq i\leq m}X_i$. Similarly, the random variable $\min_{j:\ell+1\leq j\leq n}Y_j$ is a measurable function of $\min_{j:\ell+1\leq j\leq n}X_j$. Hence, by Theorem~\ref{thm: Bernoulli}, we have that
\begin{equation*}
R\left(\min_{i:1\leq i\leq m}X_i,\min_{j:\ell+1\leq j\leq n}X_j\right)\geq R\left(\min_{i:1\leq i\leq m}Y_i,\min_{j:\ell+1\leq j\leq n}Y_j\right)=R_{m,\ell}(1-e^{-\lambda}).
\end{equation*}

\bibliographystyle{alpha}
%\bibliography{BGP}

\end{document}